\documentclass[12pt,reqno]{amsart}
\usepackage{amssymb}
\usepackage{geometry}
\usepackage[latin1]{inputenc}
\usepackage[italian,english]{babel}
\usepackage{amsmath, amsfonts, amsthm}
\usepackage{graphicx}
\usepackage{pgfplots}
\usepackage{paralist}
\usepackage{mathtools, mathabx,bigints}
\usepackage{subfig}
\usepackage{comment}
\usepackage{hyperref}
\usepackage{siunitx}
\usepackage{soul}
\numberwithin{equation}{section}
\newtheorem{thm}{Theorem}[section]

\theoremstyle{definition}
\newtheorem{rem}[thm]{Remark}
\theoremstyle{remark}

\newcommand{\norm}[1]{\left\Vert#1\right\Vert}

\newcommand{\abs}[1]{\left\vert#1\right\vert}

\newcommand{\EE}{\mathbf{E}}
\newcommand{\ix}{\hat{\mathbf{i}}_x}
\newcommand{\iy}{\hat{\mathbf{i}}_y}
\newcommand{\iit}{\hat{\mathbf{i}}_t}
\newcommand{\iid}{\hat{\mathbf{i}}_d}
\newcommand{\iin}{\hat{\mathbf{i}}_n}

\newcommand{\iz}{\hat{\mathbf{i}}_z}
\newcommand{\pp}{\mathbf{p}}
\newcommand{\xx}{\mathbf{x}}
\newcommand{\yy}{\mathbf{y}}
\newcommand{\rr}{\mathbf{r}}

\renewcommand{\Re}[1]{\operatorname{Re}\left\{#1\right\}}

\allowdisplaybreaks[3]

\let\oldoint\oint
\renewcommand{\oint}{\oldoint\nolimits}

{\left\{\begin{array}{@{}l@{}}}{\end{array}\right.}
\patchcmd{\abstract}{\scshape\abstractname}{\textbf{\abstractname}}{}{}
\makeatletter 
\def\@makefnmark{} 
\makeatother 
 
\pgfplotsset{compat=1.18}
\begin{document}
\title[\,]{Adjoint-Based Gradient Evaluation for Metasurface Inverse Design via Affine Geometric Transformations}
\author[A. Tamburrino, L. Faella, C. Forestiere, V. Mottola]{
Vincenzo Mottola$^1$, Luisa Faella$^1$, Carlo Forestiere$^2$, \\ Antonello Tamburrino$^{1}$}\footnote{\\$^1$Department of Electrical and Information Engineering, University of Cassino and Southern Lazio, Via Gaetano di Biasio 43, 03043 Cassino, Italy.\\
$^2$Department of Electrical Engineering and Information Technology,
Universit\'a degli Studi di Napoli Federico II, 80138 Napoli, Italy \\
Email: vincenzo.mottola@unicas.it {\it (corresponding author)},  antonello.tamburrino@unicas.it  {\it (corresponding author)}}

\begin{abstract} 
The sharp increasing in fabrication capabilities of nanomaterials, and complex structures such as meta-surfaces and metalens, has opened to the possibility of employing them for accurately control the electromagnetic field, beyond the possibility ensured by traditional devices. 

The demand for large scale structures and more complex functionalities from meta-surfaces lead to the research for advanced techniques of inverse design, able to conjugate the ability to produce effective designs and limited computational cost. Among the various approaches for inverse design of large meta-surfaces, the ones based on the adjoint variable method are appealing since able to ensure a minimal computational cost for the gradient computation of the cost function.

In this work, a systematic methodology for the application of the adjoint variable method for large meta-surface design is presented. The method is based on: (i) a parametrization of the relevant geometric parameters of the meta-atoms, (ii) the fast computation of the gradient with respect such parameters, allowing for the implementation of general affine transformations during the optimization process.

The main findings are first theoretically justified and a numerical validation is provided to show the effectiveness of the proposed approach. 

\end{abstract}

\maketitle
\markright{LARGE META-SURFACES DESIGN}

\section{Introduction} 
A meta-surface is a collection of a generally large number of interacting scatterers. The scatterers have dimensions comparable with the wavelength, or even smaller, while the entire structure can reach dimensions much larger of the wavelength. This inherent multiscale nature has posed significant challenges for their electromagnetic modeling and a great effort has been spent for the establishment of fast and accurate electromagnetic solvers able to compute the electromagnetic response of large meta-surfaces in a reasonable time. In this sense, methods based on fast multiple method \cite{fmm1,fmm2}, coupled mode theory \cite{cc1, cc2}, or neural networks \cite{nn1,nn2} have been successfully applied for the simulation of large area meta-surfaces.

The possibility of fast and accurate simulation of the response of a meta-surface, together with the enhanced capabilities in composites fabrication, calls for methods capable of designing this kind of structure in order to fully exploit their potential. In fact, although they have been shown to be promising in a variety of different applications such as beam steering \cite{app1,app2}, light focusing \cite{app3,app4}, bio-imaging systems \cite{app5,app6}, analog computing \cite{app7}, and imaging \cite{app8}, the functionality of the designed meta-surfaces are often limited by common problems, i.e. limited bandwidth, angular acceptance and numerical aperture. This is due the difficulty in solving electromagnetic inverse problems where a large number of unknowns are involved.

The \lq\lq classical\rq\rq approach to the design of meta-surfaces is the \lq\lq unit cell\rq\rq approach \cite{inv1,inv2}. Starting from a library of electromagnetic responses of standard meta-atoms such as pillars, disks and spheres, as a function of the designed parameters, the electromagnetic response of the meta-surface is constructed by superimposing the ones of the single \lq\lq unit cells\rq\rq. This approach requires a minimum computational cost but, at the same time, it limits the search space, requiring slowly changing structures to fulfill the underlying assumptions. Furthermore, the \lq\lq unit-cell\rq\rq approach has shown to be incompatible with the high numerical aperture requirement \cite{inv3}.

To overcame these limitations, different inverse design techniques have been proposed in literature. For example, heuristic optimization methods based on random walks across the optimization space (as particle swarm optimization or genetic algorithms) have been applied in small scale problems, where they help to discover high-performances designs \cite{heu1,heu2}. Despite this, they have been shown to be inefficient when the parameter space becomes larger \cite{heu3} and, for this reason, they are often employed in conjunction with the \lq\lq unit-cell \rq\rq approach to improve the performance of the latter \cite{heu4}.

Alternatively, gradient-based optimization scheme have been employed, due to the possibility of being integrated with an adjoint variable method for a fast gradient computation. The adjoint variable method has been already applied in different fields as fundamental tool to solve optimization problems with large number of unknowns \cite{adj1,adj2,adj3} and, recently, its usefulness as been also shown in the context of inverse design for meta-surfaces \cite{inv3,adj4,adj5}. In gradient-based optimization a certain cost functional, which depends on the design parameters, is iteratively minimized or maximized. In the iterative step, the design parameters are updated starting from the knowledge of the gradient direction. When the number of design parameters is very large (due to the huge number of meta-atoms, for example), classical finite-difference approaches for the computation of the gradient are unfeasible, since they call for the solution of a huge number of direct problems. The adjoint variable method recast the computation of the gradient as the solution of one auxiliary direct problem, where the meta-surface is illuminated by a proper excitation. In this way, for each iterative step only two different direct problems have to be solved. Such kind of approach is mandatory in meta-surfaces design due to the computational burden required by the solution of the direct problem. 

The adjoint variable method has been applied in optimization of freeform structures, where a continuous profile of dielectric permittivity is updated in order to optimize the cost function. This continuous profile is properly binarized during the process to reach feasible results for fabrication \cite{inv3,adj6,adj7,adj8}. Freeform optimization shows excellent performances but, at the same time, can pose some challenges in scaling towards large meta-surfaces. Indeed, the binarization process must be repeated several times, increasing the computational complexity. 

An alternative approach has been proposed to reduce the computational cost, based on parametrization of the single meta-atoms in terms of relevant geometrical parameters. In \cite{adj5}, for example, the meta-surface is constituted by rectangular meta-atoms, where only their width is optimized, allowing for the design of a large area meta-surface.

Another strategy for the design is represented by the integration of neural networks in the design stage. Approaches entirely based on deep learning can be found in \cite{nninv1,nninv2,nninv3}. Instead, in \cite{nninv4,nninv5}, hybrid approaches are proposed where a neural network is employed to provide a good initial design to be refined by iterative optimization scheme and adjoint variable methods.

This paper contributes in extending the capabilities of the adjoint variable method when applied to the design of large area meta-surface. In line with \cite{adj5}, the optimization variables are defined as geometrical parameters controlling the shape of each meta-atom. However, a general framework is established that can simultaneously handle multiple geometrical parameters for a single meta-atom, such as orientation, position, and size, which are varied throughout the optimization process. As a result, general geometrical transformations can be considered, since each meta-atom is subjected to a full affine transformation at every optimization step, allowing for translation, rotation, and deformation.
On the one hand, this parametrization significantly reduces the computational complexity compared to freeform optimization, as it limits the number of unknowns; on the other hand, it preserves the ability to explore a wide range of designs, owing to the generality of the geometrical transformations considered.

In Section \ref{sec:key}, the key finding of this work are summarized and the  gradient computation with respect to affine transformations is detailed. In Section \ref{sec:dim}, a rigorous justification for the results of Section \ref{sec:key} is provided. In Section \ref{sec:num}, some example of design are presented, obtained in a numerical environment. Finally, Conclusion follows in Section \ref{sec:conclusions}.

\section{Key results overview}\label{sec:key}
\begin{figure}[htb]
    \centering
    \includegraphics[width=0.5\linewidth]{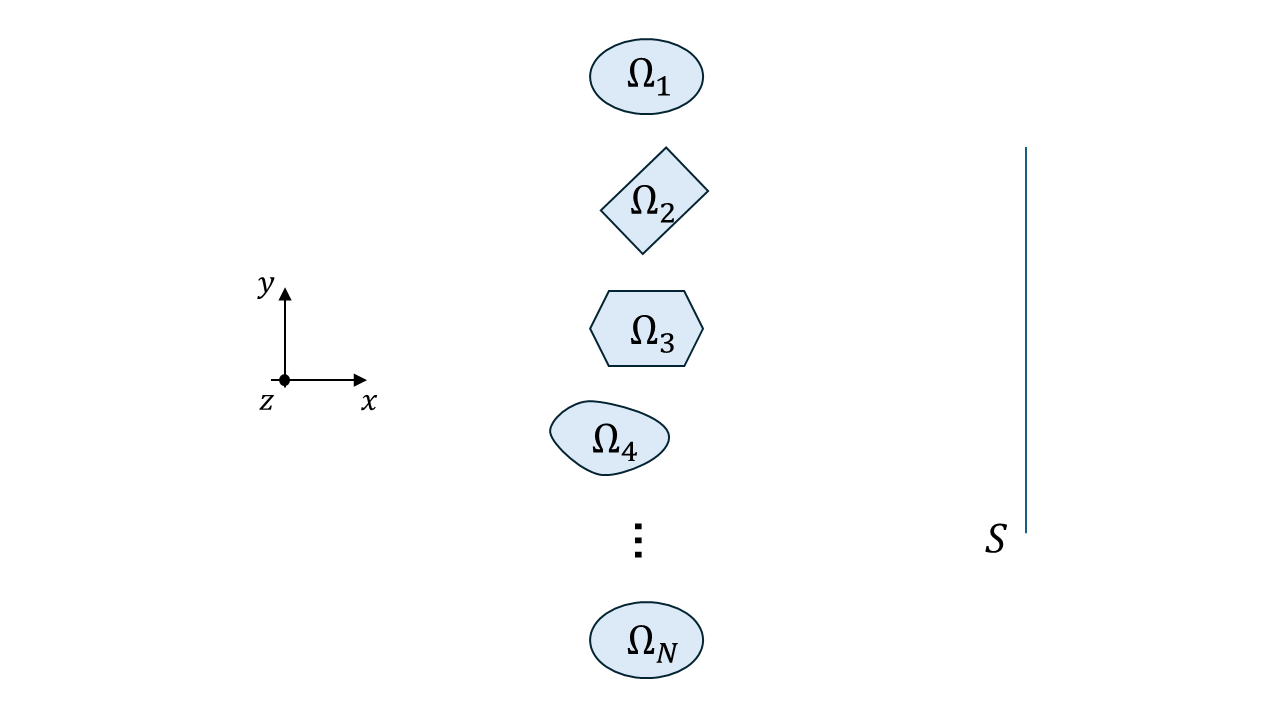}
    \caption{The meta-lens considered in the TEz propagation framework, along with the exit plane $S$ on which the target field is assigned.}
    \label{fig:tez_geo}
\end{figure}

In this paper, a TEz propagation problem is considered. The meta-surface is constituted by $N$ cylindrical scatterers of arbitrary cross-section (see Figure \ref{fig:tez_geo}). The cross section of the $i$-th meta-atom in the $xy-$plane  is denoted by $\Omega_i$, while $\varepsilon_i$ denotes its dielectric permittivity. In the following, $\EE^i(x,y)=E_x^i(x,y)\ix+E_y^i(x,y)\iy$ denotes the incident field. Since the geometry, the excitation field and the material properties are invariant in the $z$-direction, the problem is studied as a 2D problem in the $xy-$plane. 

The aim is to optimize the position and cross section shape of the each scatterer, in such a way the scattered or the total field, produced by the entire meta-surface, matches some prescribed field distribution assigned on an exit plane $S$, at a certain distance from the meta-surface. Hence, the targeted problem is an inverse design problem, in which the best scatterers' configuration is sought to achieve a prescribed target.

In doing so, iterative minimization algorithms and, in particular, gradient-based optimization scheme, are usually employed. The inverse design problem is cast in terms of minimization of a suitable cost function $I$, representing the discrepancy between the field produced by the meta-surface $\EE_f$ and the target field $\EE_d$. At each iteration, the algorithm requires the gradient of the cost function $\nabla_{\pp} I$ with respect to the design parameters $\pp$, to systematically update the solution towards a local minimum.

The core of this contribution is in a fast and accurate computation of the gradient $\nabla_{\pp} I$. Indeed, for a large area meta-surface, its computation by classical finite difference scheme would be prohibitive, since a dedicated direct simulation would be necessary for each design parameter. Hence, the attention is moved to the so-called adjoint variable method, which allows to compute the gradient, with respect to \emph{all} design parameters, by just one direct simulation.

A preliminary but fundamental aspect of the inverse design process is the definition of the design space. In this contribution, rather than adopting a freeform, or topology, optimization strategy, where the design variable would be the point-wise value of the dielectric permittivity \cite{freeform}, a parametric model is introduced, wherein the meta-atoms are characterized by a set of geometrically relevant parameters. This choice significantly reduces the dimensionality of the optimization problem compared with freeform approaches and leads to a relevant computational savings, since no binarization or filtering steps are required to obtain fabrication-ready designs. This parameterization in terms of a few relevant geometrical parameters also implicitly enforces fabrication constraints. In fact, by evolving the structure through a limited set of well-defined geometric parameters, optimization is naturally confined to realistic designs, preventing the generation of irregular or highly variable patterns that would be difficult or costly to manufacture.

The concept of parameterization was first introduced in \cite{adj5}, where rectangular meta-atoms are considered and the transformation applied during the optimization process is an isotropic expansion, i.e. the side length of the square cross section of each meta-atom was varied at each iteration.
The present work extends this approach by allowing arbitrary cross-sectional shapes and by enabling more general geometric transformations during the optimization process, thereby broadening the applicability of the method, while preserving the computation savings with respect to freeform optimization.

\begin{figure}[htb]
    \centering
    \includegraphics[width=0.75\linewidth]{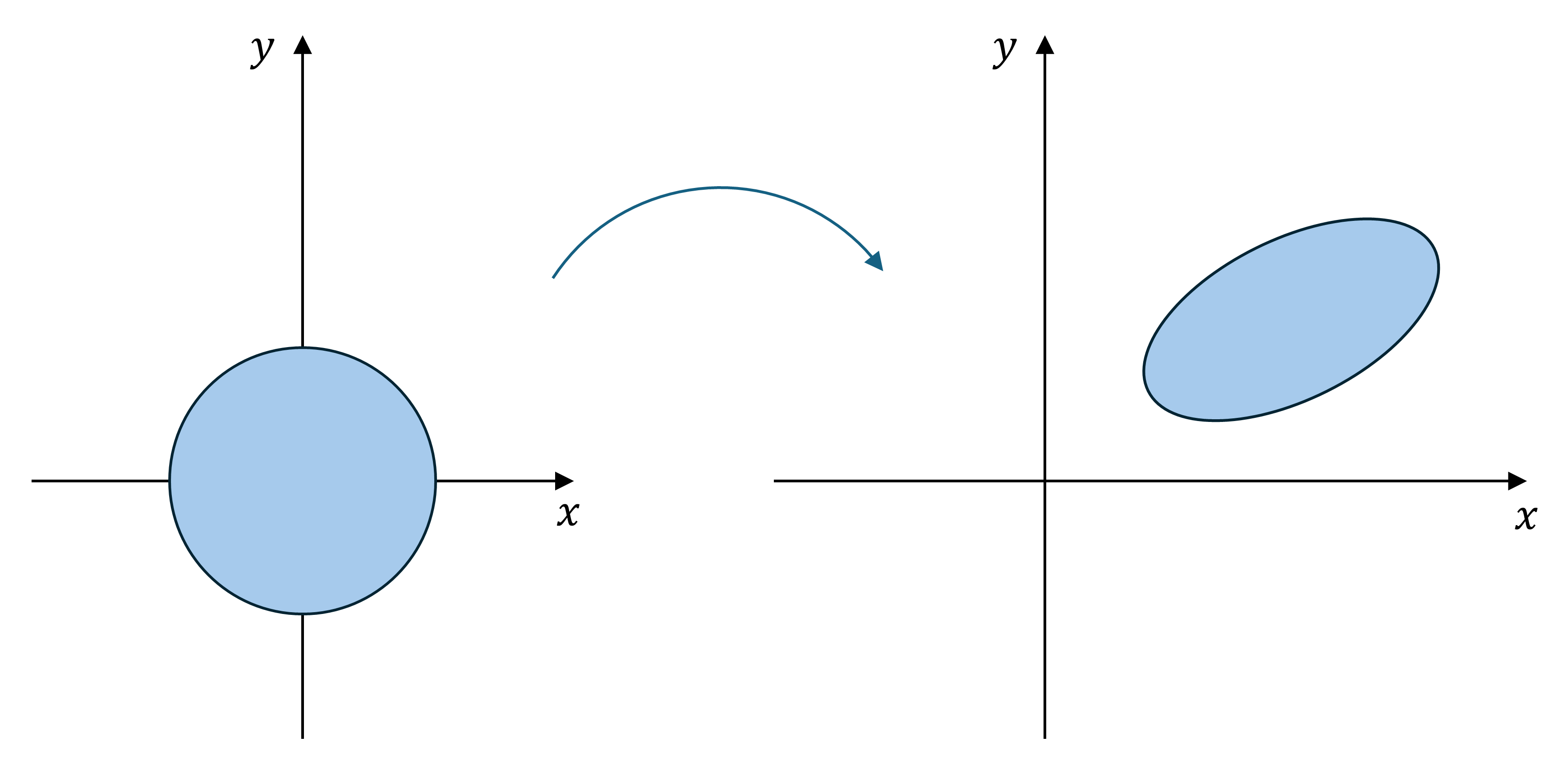}
    \caption{A generic affine transformation is applied to a circular meta-atom.}
    \label{fig:af}
\end{figure}

To be more specific, in this contribution, the shape of meta-atoms evolves under the effect of affine transformations during the optimization process, resulting from the combined action of an anisotropic deformation, a rotation, and a rigid translation (see Figure \ref{fig:af}). 
Let $\mathbf{x}\in\Omega_i$, the affine mapping applied to the $i$-th meta-atom can be expressed as
\begin{equation}
    \xx'=\begin{bmatrix}
        \cos(\theta) & \sin(\theta) \\
        -\sin(\theta) & \cos(\theta)
    \end{bmatrix}
    \begin{bmatrix}
        \lambda_x & 0 \\
        0 & \lambda_y
    \end{bmatrix}
    \mathbf{x}+\mathbf{c},
\end{equation}
where $\xx'$ are the coordinates of the transformed point, $\theta$ is the angle of the applied rotation, $\lambda_x$ and $\lambda_y$ are the anisotropic scaling factors along the two principal axes, and $\mathbf{c}$ stands for the rigid translation vector. Accordingly, each meta-atom is fully described during the optimization process by the parameter vector
\begin{equation*}
    \pp^i=[\theta^i,\lambda_{x}^i,\lambda_{y}^i,x^i,y^i],
\end{equation*}
where $\theta^i$ is the rotation angle of the $i-$th meta-atom with respect to an axis parallel to the $z$-axis and passing through its centroid, $\lambda_{x}^i$ and $\lambda_{y}^i$ control the expansion or contraction along the two references axes, and $(x^i,y^i)$ are the coordinates of the meta-atom center.

Let $p_j^i$ denote the $j-$th geometric parameter for the $i-$th meta-atom, the application of the adjoint variable method, as shown in Section \ref{sec:dim}, leads to a general expression for the computation of the gradient
\begin{equation}\label{eq_main_res}
    \frac{\partial I}{\partial p_j^i}=2\Re{\beta\oint_{\partial\Omega_i} \EE_a(\xx)\cdot \EE_f(\xx)w(\xx)\,d\xx} 
\end{equation}
where $\EE_a$ is the so-called adjoint field, while $w$ is a purely geometric weighting term, completely determined by the transformation applied to the meta-atom (rotation, anisotropic expansion, translation). The adjoint field $\EE_a$ is the response of the meta-surface when excited by a suitable adjoint source $\mathbf{J}_a$, defined only on the exit plane $S$. Both the adjoint source $\mathbf{J}_a$, and the scalar coefficient $\beta$ depend solely on the choice of the cost function.

In Equation \ref{eq_main_res}, both the adjoint field $\EE_a$ and the forward field $\EE_f$ are computed on the same geometry, that is the current design at a given iteration of the optimization process. Hence, they represent the response of the same meta-surface to different excitations. This implies that: (i) only one additional forward problem needs to be solved to compute the gradient with respect to \emph{all} the geometric parameters, for \emph{all} meta-atoms, and (ii) the computational cost of evaluating $\EE_a$ is lower than that required for $\EE_f$, since assembling the linear system underlying the numerical formulation is not necessary. 

Furthermore, Equation \ref{eq_main_res} decouples the gradient computation into two distinct problems, that are (i) determining the appropriate adjoint source, where $\mathbf{J}_a$ depends on the cost function $I$ but not on the specific transformation applied, and (ii) determining the weighting function $w$ which depends on the geometric transformation but is independent of the cost function. This decoupling has the key consequence that different geometric transformations do not require the solution of additional forward problems, but only the evaluation of \eqref{eq_main_res} with different $w$'s. Moreover, it enables the efficient reuse of the transformation weights for different objective functions.

In Section~\ref{sec:dim}, the explicit expressions of the weighting function $w$ and the adjoint source 
$\mathbf{J}_a$ are derived for the three classes of geometric transformations of interest and for a representative set of cost functions.

\section{Fast gradient computation} \label{sec:dim}
This Section is devoted to a rigorous derivation of Equation \ref{eq_main_res}. Starting from the Lorentz reciprocity theorem, the adjoint problem is formulated, and the gradient expression is recast into an appropriate integral form. This general formulation is then specialized for the three geometrical transformations of interest, and the corresponding weighting functions are analytically derived. Finally, the computation of the adjoint source is discussed for a representative set of relevant cost functions in meta-surfaces design.

\subsection{The adjoint problem}\label{sec_adj_prob}
In inverse design, the goal is to optimize the geometrical parameters of the meta-atoms so as to minimize the discrepancy between the electromagnetic field produced by the current design $\EE_f$ and a target field distribution $\EE_d$, on the exit plane $S$ of the meta-surface.

To this end, a suitable cost function $I$ is introduced to quantify the mismatch between $\EE_f$ and $\EE_d$. Since the aim is to derive a general expression for the gradient that is valid for different definitions of $I$, the adjoint problem is first formulated with respect to a fundamental metric, termed similarity function, from which the gradient for more general cost functions can be readily obtained.

Specifically, for a fixed incident field $\EE^i$, the similarity function is defined as
\begin{equation}\label{eqn:f}
    M_{\mathbf{w}}(\pp)=\int_S \mathbf{w}^*\cdot\mathbf{E}_f(\pp)\,dl,
\end{equation}
where $\mathbf{w}$ is a generic vector field defined on $S$, the superscript $^*$ denotes the complex conjugate, and $\mathbf{E}_f(\pp)$ is the total field produced by the meta-surface on the exit plane $S$ for a given set of design parameters $\pp$. The quantity $M_{\mathbf{w}}$, defined in \eqref{eqn:f}, is the scalar product between the vector $\mathbf{w}$ and the field $\mathbf{E}_f$ on $S$, and thus measures the projection of $\EE_f$ along the direction specified by $\mathbf{w}$. 

Given the definition in \eqref{eqn:f}, the goal is to derive an analytical expression for $\nabla_{\pp} M_{\mathbf{w}}$, i.e the gradient of the similarity function with respect to the design parameters. 

\begin{figure}[htb]
    \centering
    \includegraphics[width=0.9\linewidth]{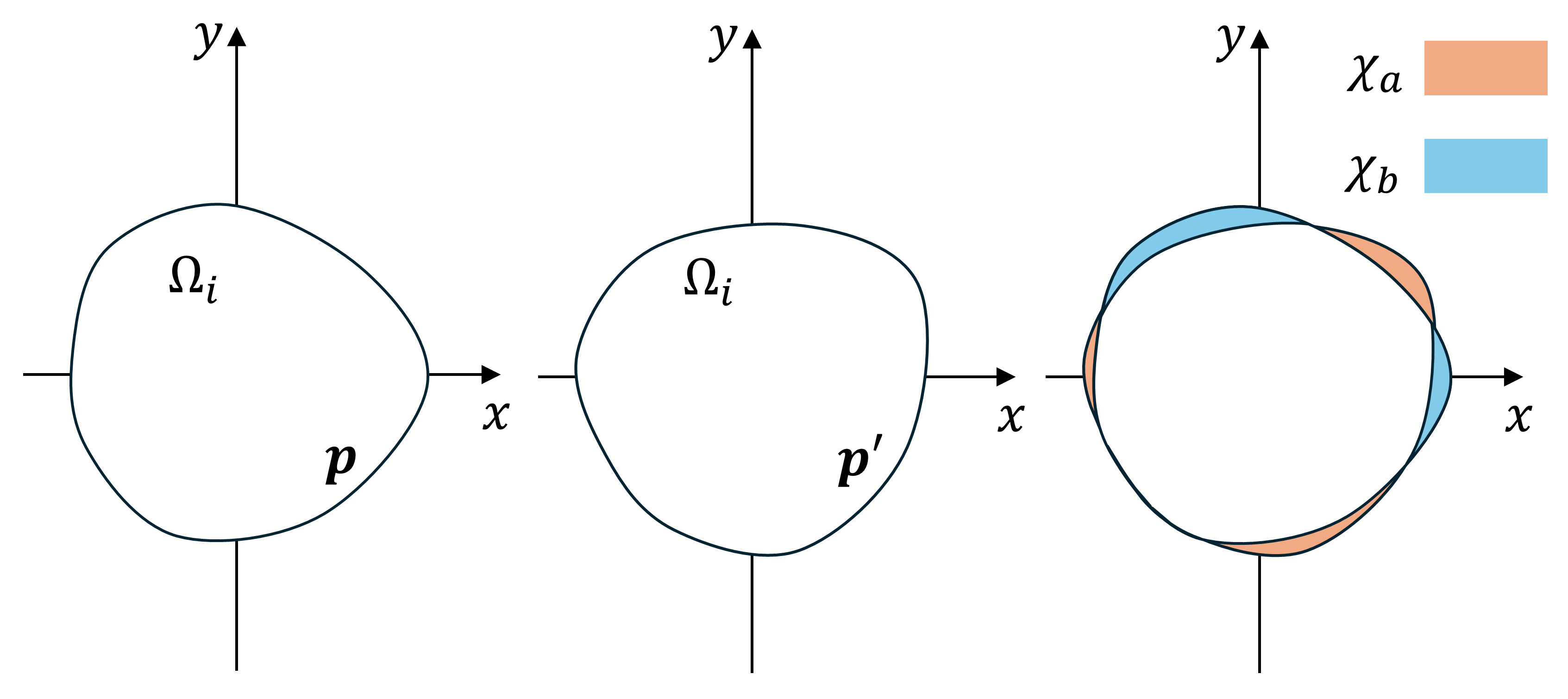}
    \caption{Left: initial design, corresponding to the vector parameter $\pp$. Center: updated design, corresponding to the vector parameter $\pp'$. Right: regions $\chi_a$ and $\chi_b$ arising from the variation in the design parameters.}
    \label{fig_setup}
\end{figure}

Consider two different sets of design parameters, namely $\pp$ and $\pp'$, as depicted in Figure \ref{fig_setup}. 
The variation $\delta M_{\mathbf{w}}$ in the similarity function due to a perturbation in the design parameters from $\pp$ to $\pp'$ is given by
\begin{equation}\label{eqn:deltai}
    \delta M_{\mathbf{w}}=M_{\mathbf{w}}(\pp')-M_{\mathbf{w}}(\pp)=\int_S \mathbf{w}^*\cdot \left(\EE_f(\pp')-\EE_f(\pp)\right)\,dl=\int_S\mathbf{w}^*\cdot\delta\EE\,dl,
\end{equation}
where $\delta\EE$ denotes the variation in the total field induced by the perturbation in the design parameters.

In what follows, $\pp'$ is obtained from $\pp$ by modifying a single meta-atom, the $i$-th one, while keeping all the others fixed. According to classic electromagnetic theory, the variation $\delta\EE$ can be expressed in terms of equivalent sources \cite{book_em}. When the geometry of the $i$-th scatterer is perturbed, two auxiliary regions can be defined (see also Figure~\ref{fig_setup})
\begin{itemize}
    \item $\chi_a^i$ is the region of the plane that belongs to the scatterer for the configuration $\pp'$ but not for $\pp$;
    \item $\chi_b^i$ is the region of the plane that belongs to the scatterer for the configuration $\pp$ but not for $\pp'$.
\end{itemize}
Then, the field variation $\delta\EE$ can be interpreted as the field radiated by an equivalent electrical current source $\mathbf{J}_{eq}$ which is different from zero only in the region $\chi^i=\chi_a^i\cup\chi_b^i$. In particular \cite{book_em},
\begin{equation}\label{eqn:source}
    \mathbf{J}_{eq}=\begin{cases}
        j\omega(\varepsilon-\varepsilon_0)\EE_f(\pp') & \text{in }\chi_a^i \\
        -j\omega(\varepsilon-\varepsilon_0)\EE_f(\pp') & \text{in }\chi_b^i \\
\mathbf{0} & \text{elsewhere},
    \end{cases}
\end{equation}
where $\EE_f(\pp')$ is the total field associated with the perturbed configuration.

The introduction of the equivalent current source $\mathbf{J}_{eq}$ provides the foundation for applying the Lorentz reciprocity theorem \cite{book_em}. This, in turn, allows the variation $\delta M_{\mathbf{w}}$ to be reformulated as an integral over regions where the meta-atom geometry changes, thereby establishing the basis for deriving the gradient expression in \eqref{eq_main_res}. Recalling that, for two source-field pairs $(\mathbf{J}_1,\EE_1)$, $(\mathbf{J}_2,\EE_2)$ defined in the same domain, the Lorentz theorem states that
\begin{equation}\label{eqn_lor_th}
        \int_{\mathbb{R}^2}\mathbf{J}_1\cdot \mathbf{E}_2 \,dS = \int_{\mathbb{R}^2}\mathbf{J}_2\cdot \mathbf{E}_1 \,dS,
\end{equation}
and by particularizing Equation \eqref{eqn_lor_th} for $\mathbf{J}_1=\mathbf{J}_{eq}$, $\EE_1=\delta \EE$, $\mathbf{J_2}=\mathbf{J}_a=\mathbf{w}^*$ and $\EE_2=\EE_a$, it follows
\begin{equation}\label{eqn:adj_prob}
\begin{split}
        \delta M_{\mathbf{w}}&=\int_S \mathbf{w}^*\cdot\delta\EE\,dl\\
        &= j\omega\left(\varepsilon-\varepsilon_0\right)\left[\int_{\chi_a^i}\EE_a(\pp)\cdot \EE_f(\pp')\,dS-\int_{\chi_b^i}\EE_a(\pp)\cdot \EE_f(\pp')\,dS\right],
\end{split}
\end{equation}
where $\EE_a$ is the total field produced by the meta-surface when the excitation is given by $\mathbf{J}_a$, defined as $\mathbf{J}_a=\mathbf{w}^*$ on the exit plane $S$ and zero elsewhere. This current is referred to as \emph{adjoint} source, and the corresponding field $\EE_a$ represents the solution of the \emph{adjoint} problem.

\begin{rem}
    Equation \eqref{eqn:adj_prob} sets the adjoint problem for the case under consideration. Moreover, it provides a basis for decoupling the effects of the chosen cost function from those of the geometrical transformation in the gradient computation. Indeed, the adjoint field $\EE_a$, and thus the adjoint source $\mathbf{J}_a$, depends only on the field $\mathbf{w}$, which is determined by the selected cost function $M_{\mathbf{w}}$. On the other hand, the geometrical transformation defines the regions $\chi_a^i$ and $\chi_b^i$, over which the integrals in \eqref{eqn:adj_prob} are evaluated.
\end{rem}

\begin{rem}
    The adjoint field $\EE_a$ is independent of the regions $\chi_a^i$ and $\chi_b^i$. Consequently, the same adjoint field can be used to compute all the components of the gradient.
\end{rem}

\subsection{The gradient of the similarity function}\label{sec_grad_comp}
In this section, the gradient of the similarity function is derived for the three different geometrical transformations considered in this paper, namely rotation, anisotropic expansion and translation. For a given meta-atom $i$, the derivative of $M_{\mathbf{w}}$ with respect to the $j-$th design parameter is evaluated as
\begin{equation}\label{eqn_der}
    \frac{\partial M_{\mathbf{w}}}{\partial p_j^i}=\lim_{\Delta p_j^i\to 0} \frac{\delta M_{\mathbf{w}}}{\Delta p_j^i}=\lim_{\Delta p_i^j\to 0}\frac{M_{\mathbf{w}}(\pp+\Delta p_{i}^i\mathbf{e}_i^j)-M_{\mathbf{w}}(\pp)}{\Delta p_i^j},  
\end{equation}
starting from the expression for $\delta M_{\mathbf{w}}$ given in \eqref{eqn:adj_prob}. In the above expression, $\mathbf{e}_i^j$ denotes the unit vector in the standard basis for the component associated with $j-$th design parameter of the $i-$th meta-atom.

As a preliminary observation, let $\Omega_i'$ be the cross section of the $i-$th meta-atom after the application of a geometrical transformation, represented by the matrix $\mathbf{T}$. Furthermore, let
\begin{equation}\label{eqn_fun_f}
    f(\xx)=j\omega(\varepsilon-\varepsilon_0)\left(\EE_a(\xx,\pp)\cdot \EE_f(\xx,\pp')\right),
\end{equation}
then the quantity $\delta M_{\mathbf{w}}$ can be conveniently rewritten as
\begin{equation}\label{eqn_der_num}
\begin{split}
    \delta M_{\mathbf{w}}=\int_{\chi} f(\xx)\,dS
    &= \int_{\Omega_i'} f(\xx)\,dS - \int_{\Omega_i} f(\xx)\,dS \\
    &=\int_{\Omega_i} \left(f(\mathbf{T}\xx)\left|\mathbf{J}_{\mathbf{T}}(\xx)\right|-f(\xx)\right)\,dS,
\end{split}
\end{equation}
where $\left|\mathbf{J}_{\mathbf{T}}\right|$ is the determinant of the Jacobian related to the transformation $\mathbf{T}$. By combing \eqref{eqn_der} and \eqref{eqn_der_num}, the gradient computation is recast to the solution of the following integral
\begin{equation}\label{eqn_der_1}
    \frac{\partial M_{\mathbf{w}}}{\partial p_j^i}=\int_{\Omega_i}\lim_{\Delta p_j^i\to 0} \frac{f(\mathbf{T}\xx)\left|\mathbf{J}_{\mathbf{T}}(\xx)\right|-f(\xx)}{\Delta p_j^i}\,dS.
\end{equation}

For the sake of clarity, the remainder of this section is organized into three subsections, each devoted to the evaluation of \eqref{eqn_der_1} with respect to a different transformation matrix $\mathbf{T}$, corresponding to the three geometric transformations considered in this work. At the end of the section, some possible generalization are presented and a more detailed discussion about function $f$ is given.

\subsubsection{Rotations}
The objective of this subsection is to compute the derivative of the similarity function with respect to a variation in the rotation angle of the $i-$th meta-atom; hence, $p_i^j=\theta^i$ in this case.
Consider the rotation of the $i-$th meta-atom about its centroid, which, without loss of generality, is assumed to coincide with the origin. The original cross section $\Omega_i$ is mapped in the perturbed geometry $\Omega_i'$, by means of the following transformation
\begin{equation*}
    \mathbf{T}=\begin{bmatrix}
        \cos(\phi) & \sin(\phi) \\ -\sin(\phi) & \cos(\phi)
    \end{bmatrix},
\end{equation*}
where $\phi$ is the angle of the applied rotation. Furthermore, for $\phi\to 0$,
\begin{equation*}
    \mathbf{T}\mathbf{x}=\begin{bmatrix}
        \cos(\phi)x + \sin(\phi)y \\ -\sin(\phi)x+\cos(\phi)y
    \end{bmatrix}=
    \begin{bmatrix}
        x + \phi y+o(\phi) \\ -\phi x+y+o(\phi)
    \end{bmatrix}.
\end{equation*}
Let $\rr=[x,y]^T$, $\hat{\mathbf{i}}_z$ be the standard unit vector in the out-of-plane direction, and neglecting higher order terms
\begin{equation}\label{eqn_rot_1}
    \mathbf{Tx}=\mathbf{x}+\phi\begin{bmatrix}
        y \\ -x
    \end{bmatrix}=
    \mathbf{x}-\phi\rr \times \hat{\mathbf{i}}_z.
\end{equation}
Plugging \eqref{eqn_rot_1} in \eqref{eqn_der_1}, and recalling that $\left|\mathbf{J}_{\mathbf{T}}\right|=1$ for this case,
\begin{equation*}
    \frac{\partial M_{\mathbf{w}}}{\partial \theta^i}
    =\int_{\Omega_i} \lim_{\phi\to 0} \frac{f(\xx-\phi\rr \times \hat{\mathbf{i}}_z)-f(\xx)}{\phi}\, dS
    =\oint_{\partial\Omega_i^-} f(\xx)\rr\cdot \iit\,dl,
\end{equation*}
where $\partial\Omega_i^-$ denotes the boundary of the $i-$th meta-atom oriented in clock-wise sense.

Therefore, the sought derivative is given by
\begin{equation}\label{eqn_der_rot_fin}
    \frac{\partial M_{\mathbf{w}}}{\partial \theta^i}=j\omega(\varepsilon-\varepsilon_0)\oint_{\partial\Omega_i^-}\EE_a\cdot\EE_f\,(\rr\cdot\iit)\,dl,
\end{equation}
which gives $w=\rr\cdot\iit$, if compared to Equation \eqref{eq_main_res}.

\subsubsection{Anisotropic expansion}
For generality, the derivative of the similarity function with respect to an expansion along an arbitrary direction is derived. Let $\hat{\mathbf{d}}=[d_x,d_y]^T$ be the unit vector defining the direction of deformation. When an expansion of modulo $\alpha$ is applied along $\hat{\mathbf{d}}$, the mapping between the original cross section $\Omega$ and the perturbed one $\Omega'$ is given by
\begin{equation}\label{eqn_exp_1}
\xx'=\mathbf{Tx}=\xx+\alpha\left(\mathbf{r}\cdot\hat{\mathbf{d}}\right)\hat{\mathbf{d}},
\end{equation}
where $\mathbf{r}=[x,y]^T$, and the matrix $\mathbf{T}$ associated with the transformation is
\begin{equation}\label{eqn_exp_2}
    \mathbf{T}=\begin{bmatrix}
        1+\alpha d_x^2 & \alpha d_x d_y \\ \alpha d_x d_y & 1+\alpha d_y^2
    \end{bmatrix}.
\end{equation}
The aim is to compute the derivative of the similarity function with respect to $p_j^i=\alpha$. By particularizing Equation \eqref{eqn_der_1}, for the transformation \eqref{eqn_exp_1}, it follows
\begin{equation*}
        \frac{\partial M_{\mathbf{w}}}{\partial \alpha}
        =\int_{\Omega_i} \lim_{\alpha\to 0} \frac{f(\xx+\alpha(\mathbf{r}\cdot\hat{\mathbf{d}})\hat{\mathbf{d}})(1+\alpha)-f(\xx)}{\alpha}\, dS
        =\oint_{\partial\Omega_i^+} f(\xx)(\rr\cdot\hat{\mathbf{d}})(\mathbf{d}\cdot \iin)\,dl.
\end{equation*}
Therefore, the sought derivative is given by
\begin{equation}\label{eqn_der_exp_fin}
    \frac{\partial M_{\mathbf{w}}}{\partial \alpha}=j\omega(\varepsilon-\varepsilon_0)\oint_{\partial\Omega_i^+}\EE_a\cdot\EE_f\,(\rr\cdot\hat{\mathbf{d}})(\mathbf{d}\cdot \iin)\,dl,
\end{equation}
which gives $w=(\rr\cdot\hat{\mathbf{d}})(\mathbf{d}\cdot \iin)$, if compared to Equation \eqref{eq_main_res}.

\subsubsection{Translation}\label{sec_tran}
In this section, the derivative of the similarity function with respect to a rigid translation of the meta-atom is considered. Specifically, the geometrical transformation is specified as follows
\begin{equation}\label{eqn_tr_1}
\xx'=\xx+\alpha\iid,
\end{equation}
where $\iid$ is a unit vector, giving the translation direction.

Following the same reasoning as for the geometrical transformations considered in the previous sections, it follows
\begin{equation*}
        \frac{\partial M_{\mathbf{w}}}{\partial \alpha}
        =\int_{\Omega_i} \lim_{\alpha\to 0} \frac{f(\xx+\alpha\iid)-f(\xx)}{\alpha}\, dS
        =\oint_{\partial\Omega_i^+} f(\xx)(\iid\cdot \iin)\,dl.
\end{equation*}
Therefore, the sought derivative is given by
\begin{equation}\label{eqn_der_tr_fin}
    \frac{\partial M_{\mathbf{w}}}{\partial \alpha}=j\omega(\varepsilon-\varepsilon_0)\oint_{\partial\Omega_i^+}\EE_a\cdot\EE_f\,(\iid\cdot \iin)\,dl,
\end{equation}
which gives $w=\iid\cdot \iin$, if compared to Equation \eqref{eq_main_res}.

\begin{rem}
   The methodology presented is fully general and can be readily extended to other geometrical transformations. 
   Furthermore, following the same approach, the method can be readily extended to three-dimensional problems as well.
\end{rem}

\subsubsection{Evaluation of the function $f$}
The aim of this section is to clarify the meaning of the function $f$, defined in \eqref{eqn_fun_f}.

The computation of the similarity function variation under perturbation of the domain is performed starting from Equation \eqref{eqn:adj_prob}, which clearly shows that $\delta M_{\mathbf{w}}$ depends on the adjoint and forward fields evaluated in the region $\chi^i=\chi_a^i\cup\chi_b^i$. 

In \eqref{eqn_fun_f}, the adjoint field $\EE_a(\xx,\pp)$ represents the response of the meta-surface to the adjoint source $\mathbf{J}_a$, when the design parameters are set to $\pp$, that is, in the unperturbed configuration. The region $\chi_a^i$ lies outside the meta-atom in this configuration; therefore, the field $\EE_a(\xx,\pp)$ is evaluated in air for $\xx \in \chi_a^i$ (see Figure~\ref{fig_setup}).
Conversely, the forward field $\EE_f(\xx,\pp')$ corresponds to the response of the meta-surface to the incident excitation when the design parameters are $\pp'$, i.e., in the perturbed configuration. Consequently, $\EE_f(\xx,\pp')$ is evaluated inside the dielectric material, characterized by $\varepsilon = \varepsilon_0 \varepsilon_r$, for $\xx \in \chi_a^i$ (see Figure~\ref{fig_setup}).  Therefore, function $f$ must be understood as
\begin{equation*}
    f(\xx)=\EE_a^{air}(\xx,\pp)\cdot\EE_f^{die}(\xx,\pp'),
\end{equation*}
where $\EE_a^{air}$ recalls that the adjoint field is evaluated in air, while $\EE_f^{die}$ stands for a field into the dielectric material. A similar reasoning applies for region $\chi_b^i$.

In Equation \eqref{eqn_der_1}, which is the basis for the gradient computation under generic geometric transformations, the integration is extended over the entire domain $\Omega_i$. The resulting expression is independent of the specific values assumed by the function $f$ in the region $\Omega \setminus \chi^i$ and can therefore be conveniently defined there.
Since the normal component of the electric field exhibits a discontinuity across the meta-atom boundary, it is necessary to ensure that the derivatives appearing in the gradient evaluation remain meaningful. For this reason, the function $f$ coincides with the physical one within the regions $\chi_a^i$ and $\chi_b^i$, that is, it takes the appropriate values of the adjoint and forward fields in air or in the dielectric, as previously discussed. Outside these regions, namely in $\Omega \setminus \chi^i$, $f$ can be smoothly and differentiably extended so as to remain continuous across the boundary $\partial\Omega$.
This choice guarantees that the gradient formulation is mathematically well-posed while preserving the correct physical field values where they are defined.

Furthermore, it should be noted that, in all the previous calculations, the partial derivative, follows from a limit of the form
\begin{multline*}
    \lim_{\alpha\to 0} \frac{f(\xx+\alpha \mathbf{v})-f(\xx)}{\alpha}\\=\lim_{\alpha\to 0} \frac{\EE_a^{air}(\xx+\alpha \mathbf{v},\pp)\EE_f^{die}(\xx+\alpha \mathbf{v},\pp')-\EE_a^{air}(\xx,\pp)\EE_f^{die}(\xx,\pp')}{\alpha},
\end{multline*}
where $\mathbf{v}$ is a proper direction dictated by the geometrical transformation considered. Since $\pp'=\pp+\alpha \mathbf{e}_i^j$, by standard arguments, it is possible to show that the above limit is equivalent to
\begin{equation*}
    \lim_{\alpha\to 0} \frac{\EE_a^{air}(\xx+\alpha \mathbf{v},\pp)\EE_f^{die}(\xx+\alpha \mathbf{v},\pp)-\EE_a^{air}(\xx,\pp)\EE_f^{die}(\xx,\pp)}{\alpha}.
\end{equation*}
On one hand this justifies the partial derivatives introduced in the calculations of the previous sections and on the other hand, it implies that all the fields appearing in the gradient computation refer to the same configuration, i.e. the unperturbed one, when the design parameters are $\pp$. 

Finally, by the interface condition for the electric field, the function $f$ appearing in \eqref{eqn_der_rot_fin}, \eqref{eqn_der_exp_fin} and \eqref{eqn_der_tr_fin} is given by
\begin{equation}\label{eqn_f_fin}
    f(\xx)=\frac{1}{\varepsilon_r}(\EE_a^{air}(\xx,\pp)\cdot\iin)(\EE_f^{air}(\xx,\pp)\cdot\iin)+(\EE_a^{air}(\xx,\pp)\cdot\iit)(\EE_f^{die}(\xx,\pp)\cdot\iit),
\end{equation}
where $\iin$, $\iit$ are the normal and tangent vector to the boundary of the meta-atom.

\subsection{The adjoint source}
In the previous section, different types of geometrical transformations of the scatterer have been addressed.
The final aspect to consider concerns the definition of the appropriate adjoint source corresponding to the desired cost function.
In the remainder of this section, several examples of cost functions are presented together with their associated adjoint sources. In what follows, $\EE_d$ denotes a target field, prescribed on the observation surface $S$, while the following notations are employed
\begin{align*}
    \langle \mathbf{u}, \mathbf{v}\rangle&=\int_{S}\mathbf{u}^*(\xx)\cdot \mathbf{v}(\xx)\,d\xx, \\
    \lVert \mathbf{u} \rVert&=\sqrt{\int_S \lvert \mathbf{u}(\xx) \rvert^2\,d\xx}, \\
    \lvert \mathbf{u}(\xx)\rvert&=\sqrt{\lvert u_x(\xx)|^2+\lvert u_y(\xx)|^2},
\end{align*}
where $u_x$ and $u_y$ are the $x$ and $y$ components of the vector $\mathbf{u}$, respectively.

\subsubsection{Magnitude of the scalar product}\label{sec_cost_1}
The considered cost function is 
\begin{equation*}
    I=|M_{\mathbf{\EE_d}}|^2,
\end{equation*}
where $M_{\mathbf{\EE_d}}$ is the similarity function defined in \eqref{eqn:f}, for $\mathbf{w}=\mathbf{E}_d$.

Preliminary, it is observed that
\begin{equation*}
    \nabla_{\pp} I= \nabla_{\pp} \left( M_{\mathbf{\EE_d}}^*M_{\mathbf{\EE_d}}\right)=2\Re{M_{\mathbf{\EE_d}}^*\nabla_{\pp} M_{\mathbf{\EE_d}}}.
\end{equation*}
Recalling the results of Section \ref{sec_adj_prob}, it follows that the computation of $\nabla_{\pp} M_{\mathbf{\EE_d}}$ requires an adjoint source given by $\mathbf{J}_a=\EE_d^*$. Therefore, the gradient is computed by applying formula \eqref{eq_main_res}, for $\mathbf{J}_a=\EE_d^*$ and $\beta=M_{\EE_d}^*$.

\subsubsection{Norm of the difference}
The considered cost function is 
\begin{equation*}
	I=\norm{\EE_d-\EE_f}^2=\langle \EE_d-\EE_f,\EE_d-\EE_f\rangle.
\end{equation*}
By rewriting the cost function as
\begin{equation*}
    \lVert \mathbf{E}_d -\mathbf{E}_f\rVert^2=\lVert \mathbf{E}_d \rVert^2 + \lVert \mathbf{E}_f \rVert^2-2\Re{\langle \mathbf{E}_d,\mathbf{E}_f\rangle},
\end{equation*}
it turns out that
\begin{equation}\label{eqn_norm_diff_tot}
    \nabla_{\pp}\left(\lVert \mathbf{E}_d -\mathbf{E}_f\rVert^2\right)=\nabla_{\pp}\lVert \mathbf{E}_f \rVert^2-2\Re{\nabla_{\pp} M_{\mathbf{E}_d}}
\end{equation}
The quantity $\nabla_{\pp} M_{\mathbf{E}_d}$ is computed accordingly to the results of Sections \ref{sec_adj_prob} and \ref{sec_cost_1}. Therefore, only the quantity $\nabla_{\pp}\lVert \mathbf{E}_f \rVert^2$ has to explicitly taken into account.

Let $G=\lVert \mathbf{E}_f(\pp)\rVert^2$, and $\delta\EE=\EE_f(\pp')-\EE_f(\pp)$, then
\begin{equation*}
    \delta G=\lVert \mathbf{E}_f(\pp') \rVert^2-\lVert \mathbf{E}_f(\pp) \rVert^2=2\Re{\langle \mathbf{E}_f(\pp),\delta\mathbf{E}\rangle}+\lVert\delta\mathbf{E}\rVert^2.
\end{equation*}
Neglecting higher order terms
\begin{equation}\label{eqn_norm_diff}
    \delta G = 2\Re{\langle \mathbf{E}_f(\pp),\delta\mathbf{E}\rangle}=2\Re{\delta M_{\EE_f}}.
\end{equation}
Equation \eqref{eqn_norm_diff} is relevant since it takes back the computation of $\nabla_{\pp} G$ to the computation of $\nabla_{\pp} M_{\EE_f}$, i.e. to the computation of the gradient of the similarity function. Hence, referring to the results of Section \ref{sec_adj_prob}, the adjoint source related to the cost function $G$ is $\mathbf{J}_a=\EE_f^*$.

Summarizing, by plugging \eqref{eqn_norm_diff} in \eqref{eqn_norm_diff_tot}
\begin{equation*}
    \nabla_{\pp}\left(\lVert \mathbf{E}_d -\mathbf{E}_f\rVert^2\right)=2\Re{\nabla_{\pp} M_{\EE_f}-\nabla_{\pp} M_{\EE_d}}=2\Re{\nabla_{\pp} M_{\EE_f-\EE_d}}.
\end{equation*}
Hence, the gradient $\nabla_{\pp} \left(\lVert \mathbf{E}_d -\mathbf{E}_f\rVert^2\right)$ is given by formula \eqref{eq_main_res}, where $\mathbf{J}_a=\EE_f^*-\EE_d^*$ and $\beta=1$.

\subsubsection{Angle between the target and forward field}
The considered cost function is
\begin{equation*}
 I=\frac{\abs{\langle \EE_d, \EE_f \rangle}^2}{\norm{\EE_d}^2\norm{\EE_f}^2}=\frac{\abs{M_{\EE_d}}^2}{\norm{\EE_d}^2\norm{\EE_f}^2}
\end{equation*}
By applying the quotient derivative rule, it holds
\begin{equation*}
\nabla_{\pp} I=\frac{\norm{\EE_f}^2\nabla_{\pp}\abs{M_{\EE_d}}^2-\abs{M_{\EE_d}}^2\nabla_{\pp}\norm{\EE_f}^2}{\norm{\EE_d}^2\norm{\EE_f}^4},
\end{equation*}
where $ \nabla_{\pp}\abs{M_{\EE_d}}^2$ is given in Section \ref{sec_cost_1}, while $\nabla_{\pp}\norm{\EE_f}^2$ is computed starting from Equation \eqref{eqn_norm_diff}. Specifically, $\nabla_{\pp}\norm{\EE_f}^2$ is given by formula \eqref{eq_main_res}, for $\mathbf{J}_a=\EE_f^*$ and $\beta=1$.

It is worth noticing that, in this case, two different adjoint problems have to be solved, one for the computation of  $ \nabla_{\pp}\abs{\langle \EE_d, \EE_f \rangle}^2$, with adjoint source $\mathbf{J}_a=\EE_d^*$ and a second one for the computation of $\nabla_{\pp}\norm{\EE_f}^2$, with adjoint source given by $\mathbf{J}_a=\EE_f^*$.

\subsubsection{Squared difference of the squared norms}
The considered cost function is
\begin{equation*}
    I=\int_{S} (\lvert \EE_d \rvert^2-\lvert \EE_f \rvert^2)^2\,dx,
\end{equation*}
The variation $\delta I$ can be computed as 
\begin{equation}\label{eqn_cost_4_1}
        \delta I=\int_S\left(\lvert \EE_f(\mathbf{p}')\rvert^4-\lvert \EE_f(\mathbf{p})\rvert^4\right)\,dx
-2\int_S \lvert \EE_d\rvert^2\left(\lvert \EE_f(\mathbf{p}')\rvert^2-\lvert \EE_f(\mathbf{p})\rvert^2\right)\,dx
\end{equation}
Let $\EE_f(\pp')=\EE_f(\pp)+\delta\EE$, then
\begin{equation}\label{eqn_cost_4_2}
           A_1=\lvert \EE_f(\mathbf{p}')\rvert^4-\lvert \EE_f(\mathbf{p})\rvert^4
            = 4\Re{\lvert\EE_f(\mathbf{p})\lvert^2 \EE_f^*(\mathbf{p})\cdot\delta\EE} + o\left(\lvert \delta\EE\rvert\right)
\end{equation}
and
\begin{equation}\label{eqn_cost_4_3}
        A_2=\lvert \EE_d\rvert^2\left(\lvert \EE_f(\mathbf{p}')\rvert^2-\lvert \EE_f(\mathbf{p})\rvert^2\right)       
        =2\Re{\lvert \EE_d\rvert^2\mathbf{E}_f^*(\mathbf{p})\cdot\delta\EE}+o\left(\lvert \delta\EE\rvert\right)
\end{equation}
Combining Equation \eqref{eqn_cost_4_1}, \eqref{eqn_cost_4_2} and \eqref{eqn_cost_4_3}, it holds
\begin{equation}\label{eqn_cost_4_4}
\delta I=4\Re{\int_S\left(\abs{\EE_f(\pp)}^2-\abs{\EE_d(\pp)}^2\right) \EE_f^*(\mathbf{p})\cdot\delta\EE\,dx}.
\end{equation}
where higher order terms have been neglected. Equation \eqref{eqn_cost_4_4} leads to
\begin{equation*}
    \delta I = 4\Re{\delta M_{\mathbf{w}}},
\end{equation*}
with $\mathbf{w}=\left(\abs{\EE_f(\pp)}^2-\abs{\EE_d(\pp)}^2\right) \EE_f$. Hence, the gradient is given again by formula \eqref{eq_main_res}, where the adjoint source is $\mathbf{J}_a=\left(\abs{\EE_f(\pp)}^2-\abs{\EE_d(\pp)}^2\right) \EE_f^*(\pp)$ and $\beta=2$.

\subsubsection{Angle between the square of the norms}
The considered cost function is
\begin{equation*}
    I=\frac{\left( \int_S \lvert \EE_d\rvert^2\lvert\EE_f\rvert^2\,dx\right)^2}{\int_S \lvert \EE_d \rvert^4\,dx\int_S \lvert \EE_f \rvert^4\,dx}.
\end{equation*}
Let
\begin{equation*}
    F=\int_S \lvert \EE_d \rvert^2 \lvert \EE_f\rvert^2\,dx,\quad\quad G=\int_S \lvert \EE_f \rvert^4\,dx, \quad\quad H=\int_S \lvert \EE_d \rvert^4\,dx
\end{equation*}
by applying the quotient rule for differentiation, it holds
\begin{equation*}
    \nabla_{\pp} I=\frac{\left(2 F \nabla_{\pp} F\right)G-F^2\nabla_{\pp} G}{G^2 H}.
\end{equation*}
The quantity $\nabla_{\pp} F$ is computed starting from Equation \eqref{eqn_cost_4_3}, which gives $\mathbf{J}_a=\lVert\EE_d\rVert^2\EE_f^*$ and $\beta=1$; while $\nabla_{\pp} G$ is computed starting by Equation \eqref{eqn_cost_4_2}, obtaining $\mathbf{J}_a=\lVert\EE_f\rVert^2\EE_f^*$ and $\beta=2$.

\subsubsection{Square norm of the field at a given point} 
The considered cost function is
\begin{equation*}
	I=\abs{\EE_f(\xx_0)}^2,
\end{equation*}
where $\xx_0$ is a given point in which the norm of field $\EE_f$ has to be maximized.

In this case, the variation $\delta I$ is given by
\begin{equation*}
\delta I = \abs{\EE_f(\pp')}^2-\abs{\EE_f(\pp)}^2=2\Re{\EE_f(\pp)\cdot\delta\EE}+\abs{\delta\EE}^2.
\end{equation*}
Noticing that
\begin{equation*}
    \delta I = 2 \Re{\delta M_{\mathbf{w}}},
\end{equation*}
for $\mathbf{w}=\mathbf{E}_f(\mathbf{x}_0)\delta(\mathbf{x}-\mathbf{x}_0)$, it follows that the sought adjoint source is given by
\begin{equation*}
    \mathbf{J}_a=\mathbf{E}_f^*(\mathbf{x}_0)\delta(\mathbf{x}-\mathbf{x}_0),
\end{equation*}
while $\beta=1$. In other words, for this special case, the adjoint source is given by a singular source, i.e. a current dipole placed in $\xx_0$.

\section{Numerical validation}\label{sec:num}
In this section, the effectiveness of the proposed approach is investigated through numerical experiments. First, a comparison between the gradient computed via the adjoint method and that obtained with a finite-difference scheme is presented, in order to validate the accuracy of the gradient computation. Then, different inverse design examples are discussed, where the proposed method is integrated into an iterative gradient-based optimization framework.

\subsection{Case of study description}\label{sec_geo_sim}
\begin{figure}[htb]
    \centering
    \includegraphics[width=0.35\linewidth]{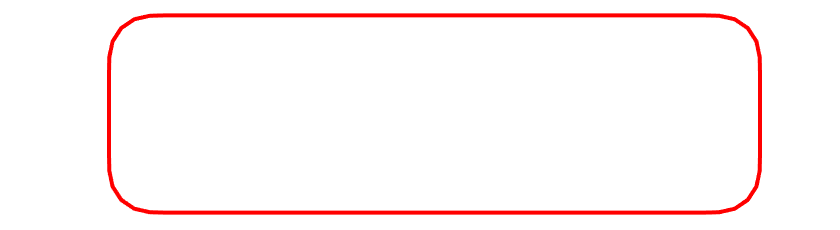}
    \caption{Geometry of the scatterers considered in the simulations.}
    \label{fig_geo_meta}
\end{figure}

Throughout this section, a specific geometry is considered for the meta-atoms forming the analyzed metasurfaces, as illustrated in Figure~\ref{fig_geo_meta}.
Each scatterer consists of a rectangular element with rounded corners, obtained by replacing each sharp corner with a circular arc. The rectangle has a width of $L_x=\qty{660}{\nano\meter}$ and a height of $L_y=\qty{200}{\nano\meter}$. The arc at the top-right corner $(L_x/2,L_y/2)$ has center at $(L_x/2-R,L_y-R)$ with radius $R=L_x/8$ and is tangent to the lines $x=L_x/2$ and $y=L_y/2$. The other corners are treated analogously.

The dielectric material filling the meta-atom is titanium dioxide (TiO$_2$), characterized by a relative permittivity of $\varepsilon_r = 5.76$ (corresponding to a refractive index $n = 2.4$) in the frequency range of interest \cite{mat01}.

The free-space wavelength is set to $\lambda_0=\qty{660}{\nano\meter}$, which corresponds to red light.

For all the simulations, the incident field is a plane wave, polarized in the $x$ direction
\begin{equation}\label{eqn_inc}
    \EE^i(\xx)=E^ie^{-jk_0x}\hat{\mathbf{i}}_y,
\end{equation}
where $k_0=2\pi/\lambda_0$.

\subsection{Numerical model}\label{sec_num}
The numerical model for solving the TEz scattering problem is derived following the boundary element method formulation presented in \cite{bem85}. 

In this formulation, the unknowns are the tangential component of the electric and magnetic fields on the boundary of each meta-atom. However, to compute the gradient via the adjoint method (see Equations~\eqref{eq_main_res} and~\eqref{eqn_f_fin}), the normal component of the electric field is also required on the meta-atom boundaries. This quantity can be readily obtained from the boundary values of the unknowns, specifically, from the magnetic field distribution on the boundary. Recalling that, in the TEz setting, $\mathbf{H}=H\iz$, it follows
\begin{equation}\label{eqn_en}
    \EE\cdot\iin=\frac{1}{j\omega\varepsilon_0}\nabla\times\mathbf{H}\cdot\iin=\frac{1}{j\omega\varepsilon_0}\nabla H\times\iz \cdot \iin=\frac{1}{j\omega\varepsilon_0}\nabla H \cdot \iit=\frac{1}{j\omega\varepsilon_0}\partial_t H,
\end{equation}
where $\iin$ is the outward normal to the meta-atom boundary, $\iit$ is the tangent unit vector counterclockwise oriented and $\partial_t$ denotes the tangent derivative. In \eqref{eqn_en}, the second equality follows from the identity $\nabla\times(\alpha\mathbf{u})=\nabla\alpha\times\mathbf{u}+\alpha\nabla\times\mathbf{u}$ and $\nabla\times\mathbf{\iz}=\mathbf{0}$, while the third equality results from a circular shift in the triple product and $\iz\times\iin=\iit$. Hence, the normal component of the electric field can be efficiently computed from the magnetic field with negligible computational cost. 

The other quantity needed for the gradient computation is the electric field $\EE_a^i$ generated by the adjoint source, that is the field produced by the electrical current density $\mathbf{J}_a$ placed on the line $S$. This is obtained starting from the knowledge of the Green's functions for the TEz problem, that are \cite{Richmond1966TE}
\begin{align*}
    G_{11}&=-\frac{k}{4\omega\varepsilon_0\rho^3}\left[k\rho y^2H_0^{(2)}(k\rho)+(x^2-y^2)H_1^{(2)}(k\rho)\right]\\
    G_{12}&=-\frac{k}{4\omega\varepsilon_0\rho^3}\left[2H_1^{(2)}(k\rho)-k\rho H_0^{(2)}(k\rho)\right]xy\\
    G_{21}&=-\frac{k}{4\omega\varepsilon_0\rho^3}\left[2H_1^{(2)}(k\rho)-k\rho H_0^{(2)}(k\rho)\right]xy\\
    G_{22}&=-\frac{k}{4\omega\varepsilon_0\rho^3}\left[k\rho x^2H_0^{(2)}(k\rho)+(y^2-x^2)H_1^{(2)}(k\rho)\right],
\end{align*}
where $\rho=\lVert \yy-\xx\rVert$, with $\xx$ observation point and $\yy$ spatial position of the source; $x$ and $y$ are the $x$ and $y$ components of the vector $\xx-\yy$. The sought field is then obtained by integrating on the line $S$
\begin{equation*}
    \EE_a^i(\xx)=\int_S \mathbf{G}(\xx,\yy)\mathbf{J}_a(\yy)\,d\yy,
\end{equation*}
where 
\begin{equation*}
    \mathbf{G}=\begin{bmatrix}
        G_{11} & G_{12} \\ G_{21} & G_{22}
    \end{bmatrix}.
\end{equation*}

\subsection{Gradient accuracy}\label{sec_grad_acc}
In this section, the gradient computed using Equation~\eqref{eq_main_res} is compared with that obtained through a finite difference scheme, in order to validate the proposed method. The analysis is carried out in the presence of three scatterers, all sharing the same geometry described in Section~\ref{sec_geo_sim}. The scatterers are aligned along the $y-$axis, with a vertical spacing of $d=\qty{726}{\nano\meter}$ between their centroids. The centroid of the central meta-atom is placed at the origin, while the remaining two are located at $(0,\pm d)$.

The choice of the target field is, instead, inspired by light focusing applications. In particular, the aim is to design a meta-surface able to focusing light in a prescribed point $\xx_0$. A possible choice for $\EE_d$ is the time-reversal of the field produced by an electrical dipole placed in $\xx_0$. Hence, let $\EE^{dip}(\xx,\xx_0)$ be the electric field generated by the dipole placed in $\xx_0$ and directed along the $y-$axis, i.e. $\EE^{dip}(\xx,\xx_0)=\mathbf{G}(\xx,\xx_0)\cdot J_0\mathbf{\hat{i}}_y$, the target field is $\EE_d=(\EE^{dip})^*|_S$. For the results of this section, $\xx_0=[37\lambda_0,\lambda_0]^T$ and $J_0=\qty{1}{\ampere\per\square\meter}$ and $S$ it the straight line from the point $(\qty{1320}{\nano\meter},\qty{-1056}{\nano\meter})$ to the point $(\qty{1320}{\nano\meter},\qty{1056}{\nano\meter})$.

The accuracy of the proposed gradient computation is assessed by varying the design parameters of the central meta-atom, whose centroid is located at the origin, while keeping the parameters for the other two fixed. To evaluate the method over a representative range of configurations, the gradient is computed for several values of the design parameters. Specifically, let $\pp=[\theta,L_x,y_c]$, as the considered design parameters, where $\theta$ is rotation angle of the meta-atom, $L_x$ is its width, and $y_c$ is the $y-$ coordinate of the centroid, the derivative of the cost function is evaluated:
(i) with respect to a rotation of the meta-atom, where the rotation angle $\theta$ is varied within the range $[0,\pi/2]$;
(ii) with respect to a deformation along the $x$-axis, corresponding to the transformation in \eqref{eqn_exp_1} with $\hat{\mathbf{d}}=\hat{\mathbf{i}}_x$, where the meta-atom width $L_x$ is varied from $\qty{330}{\nano\meter}$ to $\qty{660}{\nano\meter}$; and
(iii) with respect to a translation along the $y$-axis, corresponding to the transformation in \eqref{eqn_tr_1} with $\hat{\mathbf{i}}_d=\hat{\mathbf{i}}_y$, where the centroid position $y_c$ is shifted from $[\qty{0}{\nano\meter},\qty{-50}{\nano\meter}]$ to $[\qty{0}{\nano\meter},\qty{-50}{\nano\meter}]$. 

\begin{figure}[htb]
 \centering
\includegraphics[width=.9\textwidth]{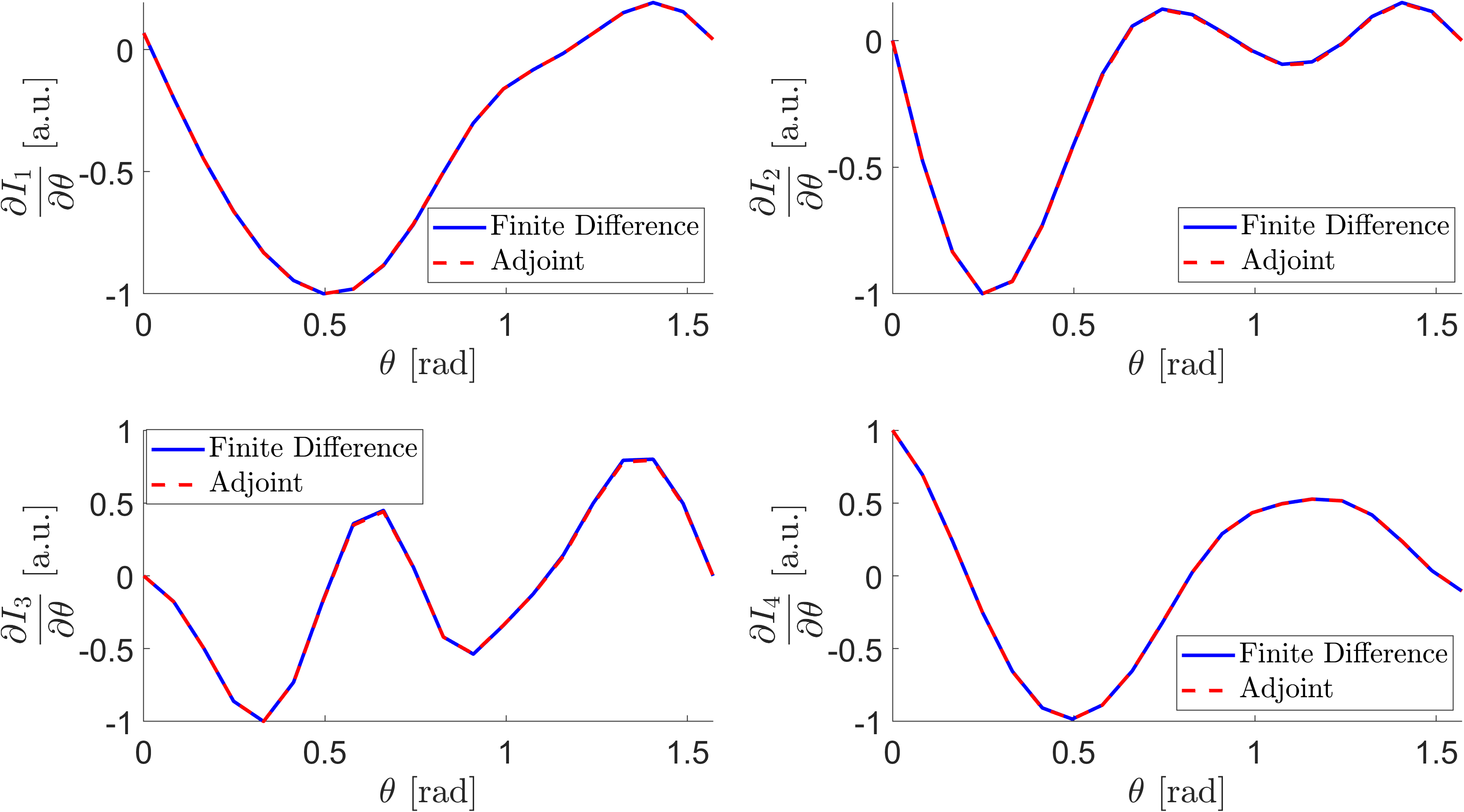} 
 \caption{Comparison between the gradient with respect to a rotation in the meta-atom, computed via adjoint method and a finite difference scheme. The gradient is computed for different values of the rotation angle.}
 \label{fig_grad_acc_1}
 \end{figure}

\begin{figure}[htb]
 \centering
\includegraphics[width=.9\textwidth]{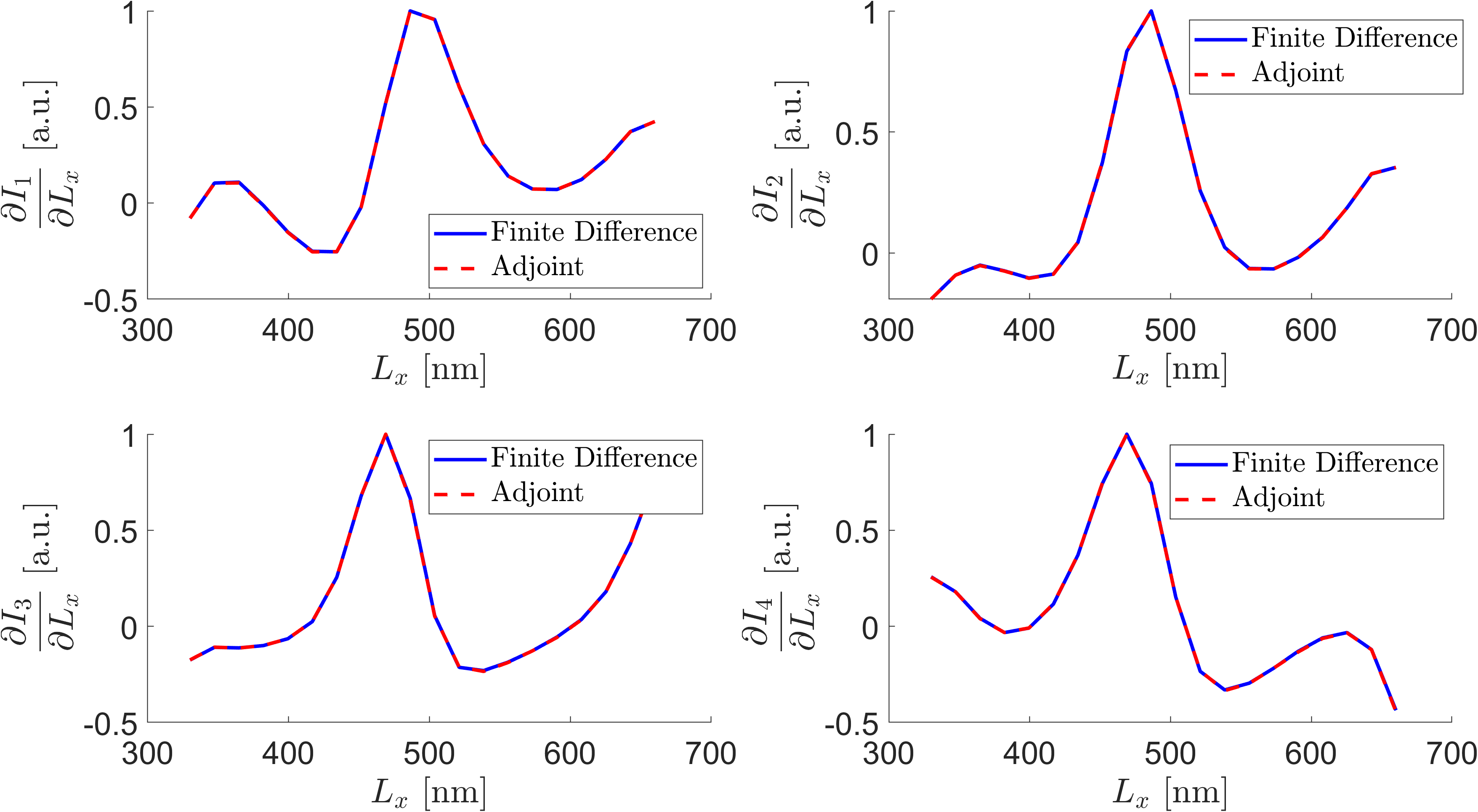} 
 \caption{Comparison between the gradient with respect to the meta-atom width, computed via adjoint method and a finite difference scheme. The gradient is computed for different values of the meta-atom width.}
 \label{fig_grad_acc_2}
 \end{figure}

\begin{figure}[htb]
 \centering
\includegraphics[width=.9\textwidth]{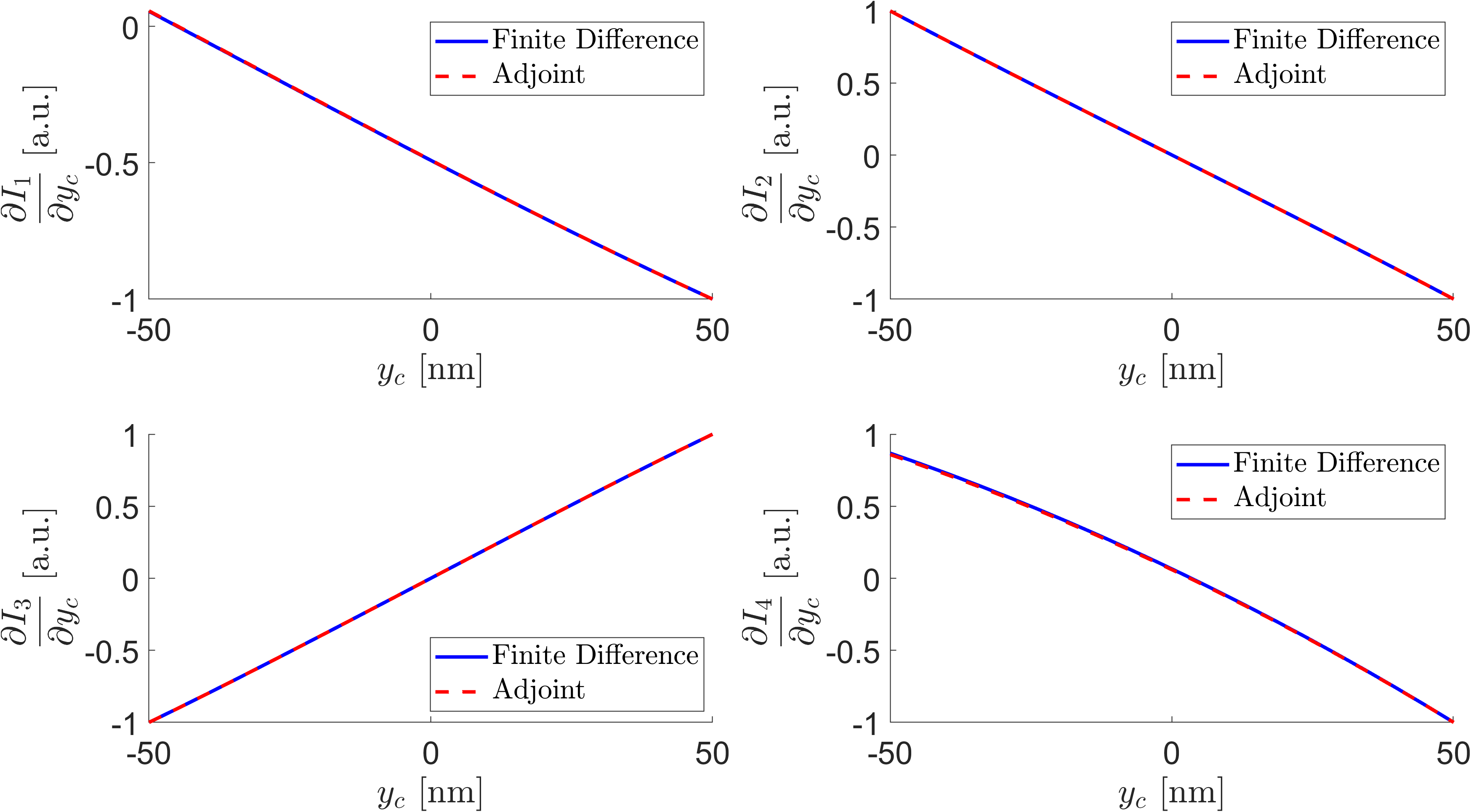} 
 \caption{Comparison between the gradient with respect to the $y$ coordinate of the centroid, computed via adjoint method and a finite difference scheme. The gradient is computed for different values of $y$ coordinate of the centroid.}
 \label{fig_grad_acc_3}
 \end{figure}

The results of the comparison are reported in Figures \ref{fig_grad_acc_1}, \ref{fig_grad_acc_2} and \ref{fig_grad_acc_3}, where different combinations of cost function and geometrical transformation are reported. Specifically, the cost functions $I_1=|M_{\EE_d}|^2$, $I_2=\lVert \EE_f - \EE_d\rVert^2$, $I_3=\lVert \lvert \EE_f\rvert^2-\lvert \EE_d\rvert^2\rVert$, and $I_4=\lvert E_f(\xx_0)\rvert^2$ are analyzed.
As it can be seen, the two different methods closely match, confirming the accuracy of the proposed approach. 

\subsection{Inverse design examples}
The fast calculation of the gradient proposed in this contribution is embedded in an iterative optimization scheme to solve some demonstrative inverse design problems. 

In the following, the adopted cost function is
\begin{equation*}
    I=\frac{\abs{\langle \EE_d,\EE_f\rangle}^2}{\norm{\EE_D}^2\norm{\EE_f}^2},
\end{equation*}
while the design parameters are assumed to be the rotation angle of each meta-atom constituting the meta-surface.
The optimization is carried out by employing the MATLAB's built-in function \texttt{fminunc} \cite{MATLAB2025}, which implements the BFGS Quasi-Newton method \cite{NocedalWright2006}.

The incident field is the one reported in \eqref{eqn_inc}, i.e. a plane wave propagating in the positive $x-$direction, while the free-space wavelength is $\lambda_0=\qty{660}{\nano\meter}$ (red light).

The analyzed meta-surfaces is constituted by 128 meta-atoms, which share the same geometry depicted in Figure \ref{fig_geo_meta}. The width of the meta-atoms is set to $L_x=\qty{660}{\nano\meter}$ and the height to $L_y=\qty{200}{\nano\meter}$. The rotation angles are initialized to zero, that is the main axis of each meta-atom is parallel to the $x-$axis. The meta-atoms are aligned vertically and the centroids are spaced of $d=\qty{726}{\nano\meter}$. The meta-atoms are filled by Ti0$_2$ oxide, with a refractive index equal to $2.4$. The whole meta-surface has a length of $\qty{92.9}{\micro\meter}$, which corresponds to $141\lambda_0$.

\begin{figure}[htb]
    \centering
    \includegraphics[width=0.85\linewidth]{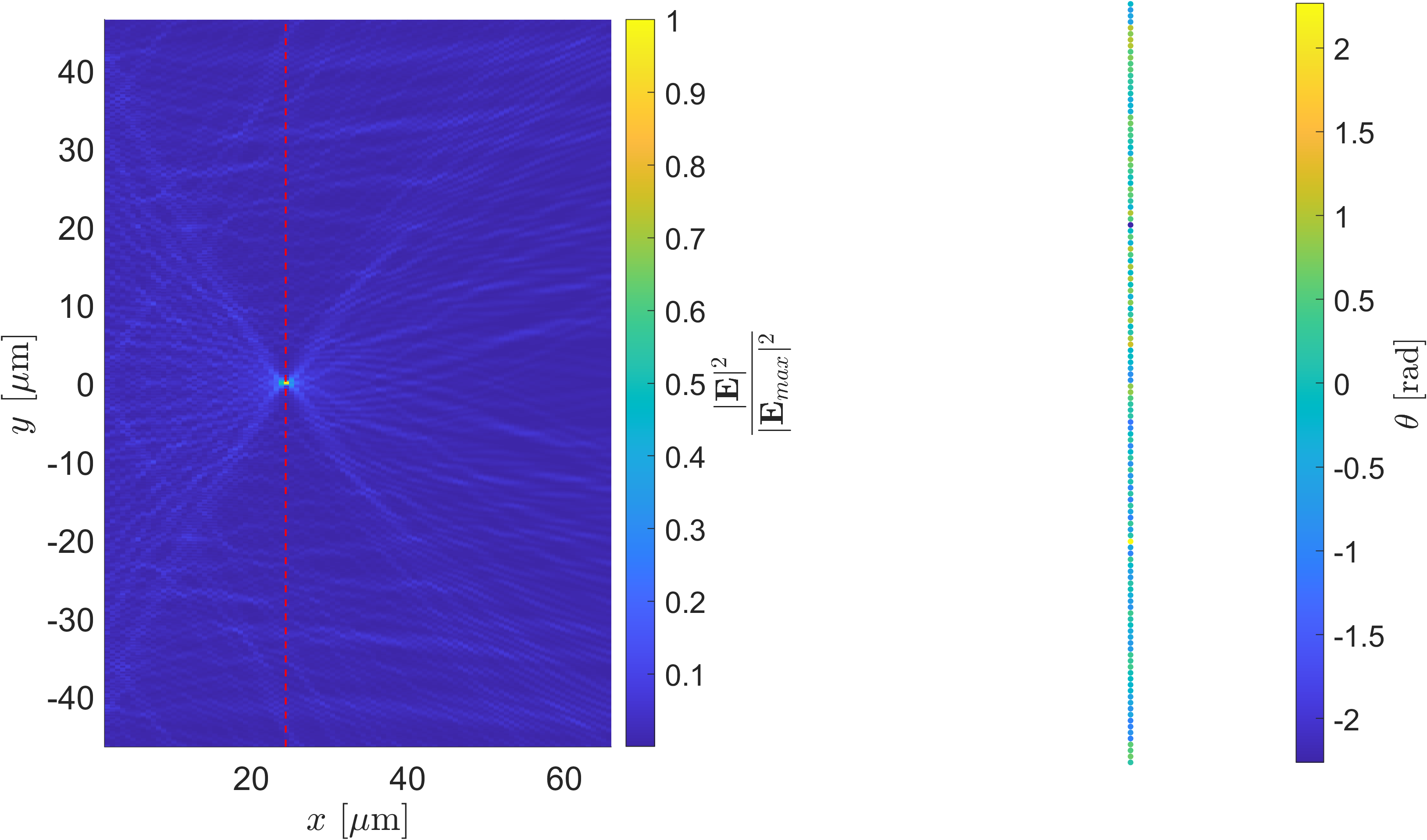}
    \caption{Left: Normalized intensity of the total electric field. The red dashed line is at $37\lambda_0$. Right: Optimal values of the rotation angle for each scatter.}
    \label{fig_res1_1}
\end{figure}

\begin{figure}[htb]
    \centering
    \includegraphics[width=0.85\linewidth]{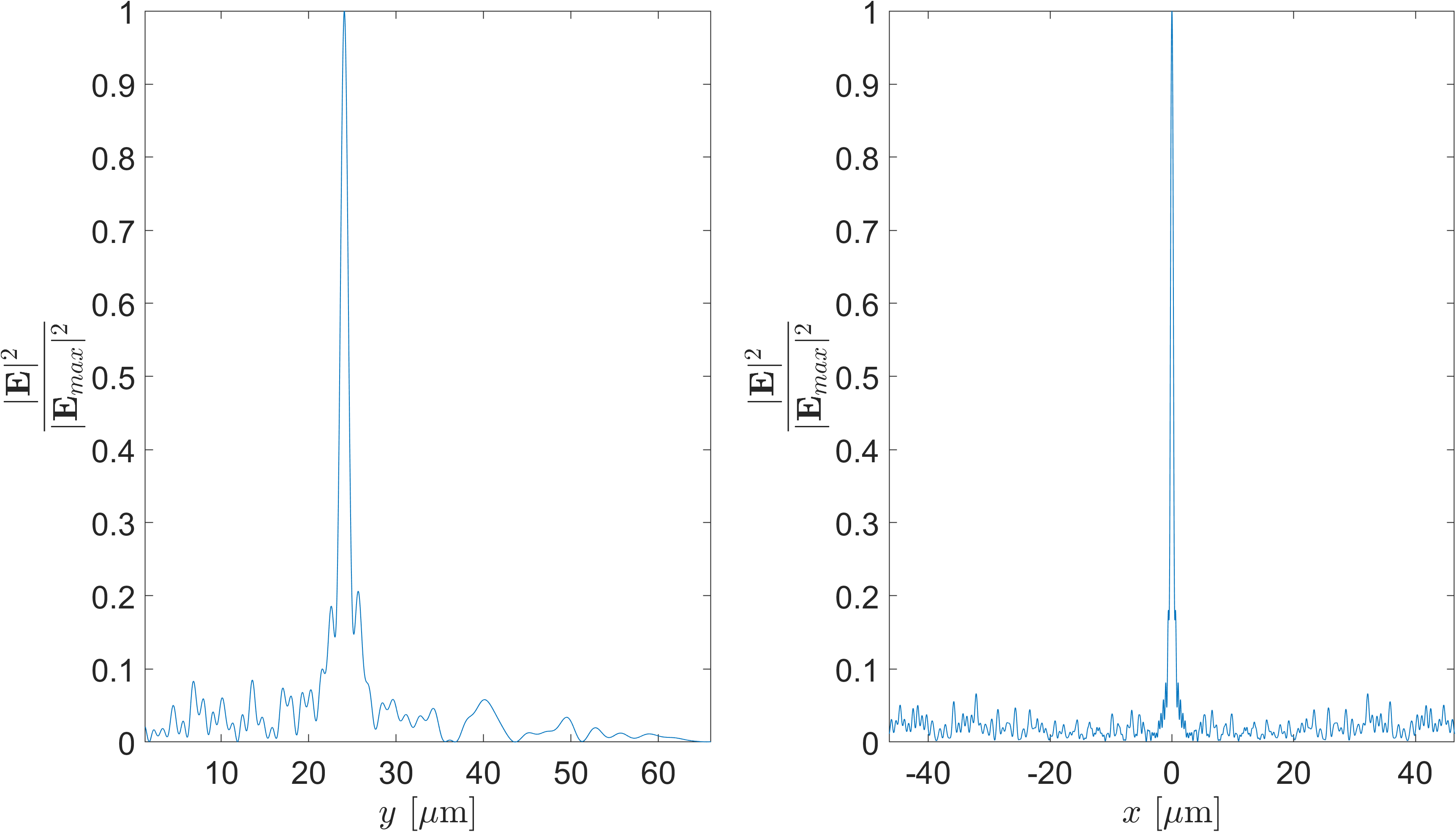}
    \caption{Left: Normalized intensity of the electric field as a function of the distance from the meta-surface. Right: Normalized intensity of the electric field at the focal spot, as a function of the vertical coordinate.}
    \label{fig_res1_2}
\end{figure}

As first inverse design example, the optimal design able to focalize light in a prescribed point is determined. The point is located at $(37\lambda_0,0)$. As already discussed in Section \ref{sec_grad_acc}, the target field distribution is given by the time-reversal of the electric field radiated by a unitary current dipole aligned along the $y-$axis and placed in the focal point. The optimization results are summarized in Figures \ref{fig_res1_1} and \ref{fig_res1_2}.
In Figure \ref{fig_res1_1}(a), the normalized electric field intensity is shown. A distinct focal spot is clearly visible, located at a distance of $36.65\lambda_0$, which closely matches the target focal length. It is worth noting that the optimal design exhibits a highly non-symmetric and rapidly varying pattern, with rotation angles that change sharply along the sequence of meta-atoms. In Figure \ref{fig_res1_2}(a), the electric field intensity is shown as a function of the $y-$coordinate at the focal plane, while Figure \ref{fig_res1_2}(b) presents the electric field intensity as a function of the $x-$coordinate for $y=0$. The full width of the main lob at half maximum (FWHM) is $0.84\lambda_0=\qty{551}{\nano\meter}$.

\begin{figure}[htb]
    \centering
    \includegraphics[width=0.85\linewidth]{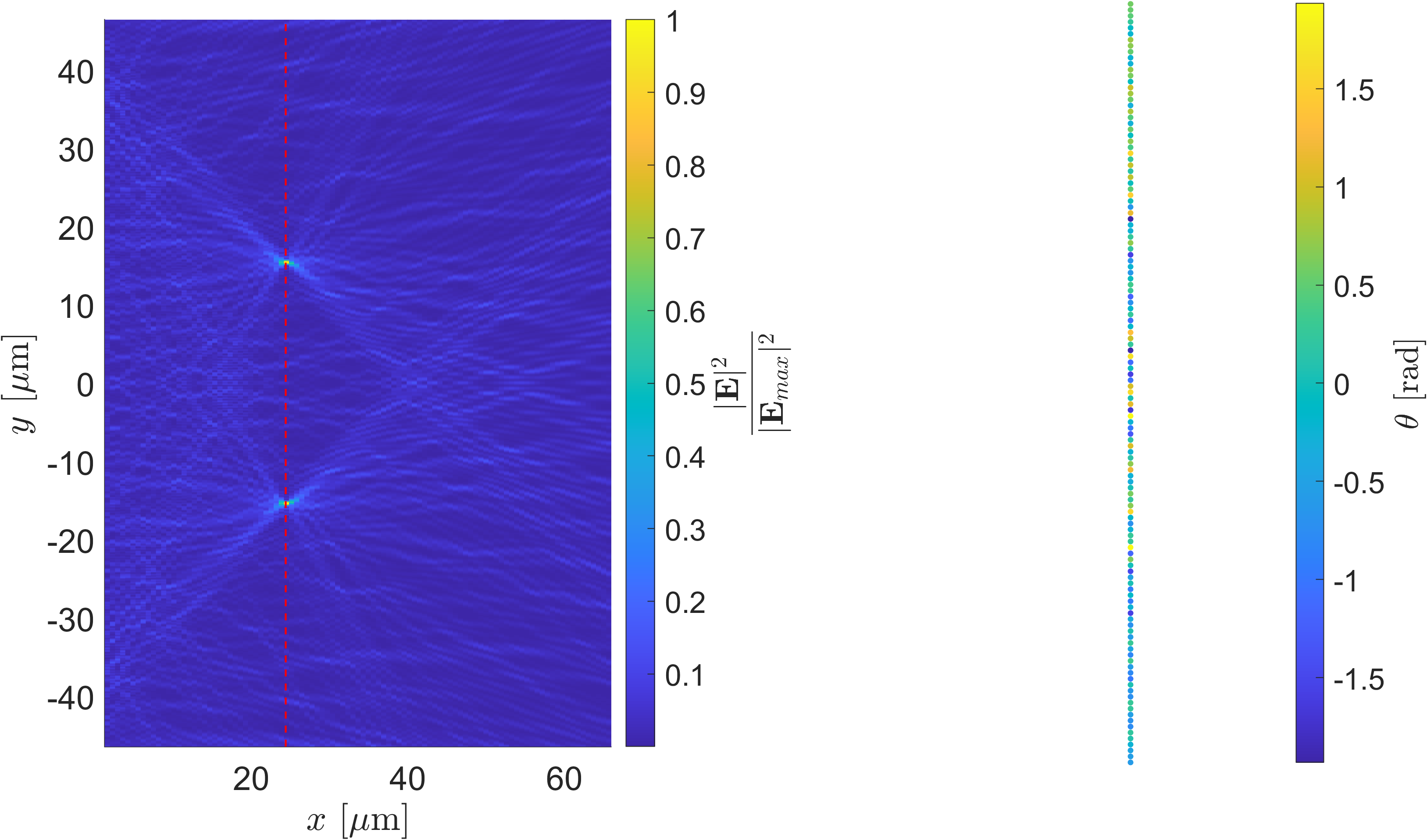}
    \caption{Light focusing in three focal spots. Left: Normalized intensity of the total electric field. The red dashed line is at $37\lambda_0$. Right: Optimal values of the rotation angle for each scatter.}
    \label{fig_res2_1}
\end{figure}

To further assess the feasibility of the proposed method, additional optimizations were performed to focus light simultaneously at multiple locations.
Figure \ref{fig_res2_1} shows the results obtained for two prescribed focal points. The target field distribution in this case is given by
\begin{equation}\label{eqn_super_dip}
    \EE_d(\xx)=\left(\EE^{dip}(\xx,\xx_0)+\EE^{dip}(\xx,\xx_1)\right)^*, \quad \text{for }\xx\in S,
\end{equation}
where $\xx_0=(37\lambda_0,23.45\lambda_0)$, $\xx_1=(37\lambda_0,-23.45\lambda_0)$ and $S$ is the same exit plane of the previous example. 

\begin{figure}[htb]
    \centering
    \includegraphics[width=0.85\linewidth]{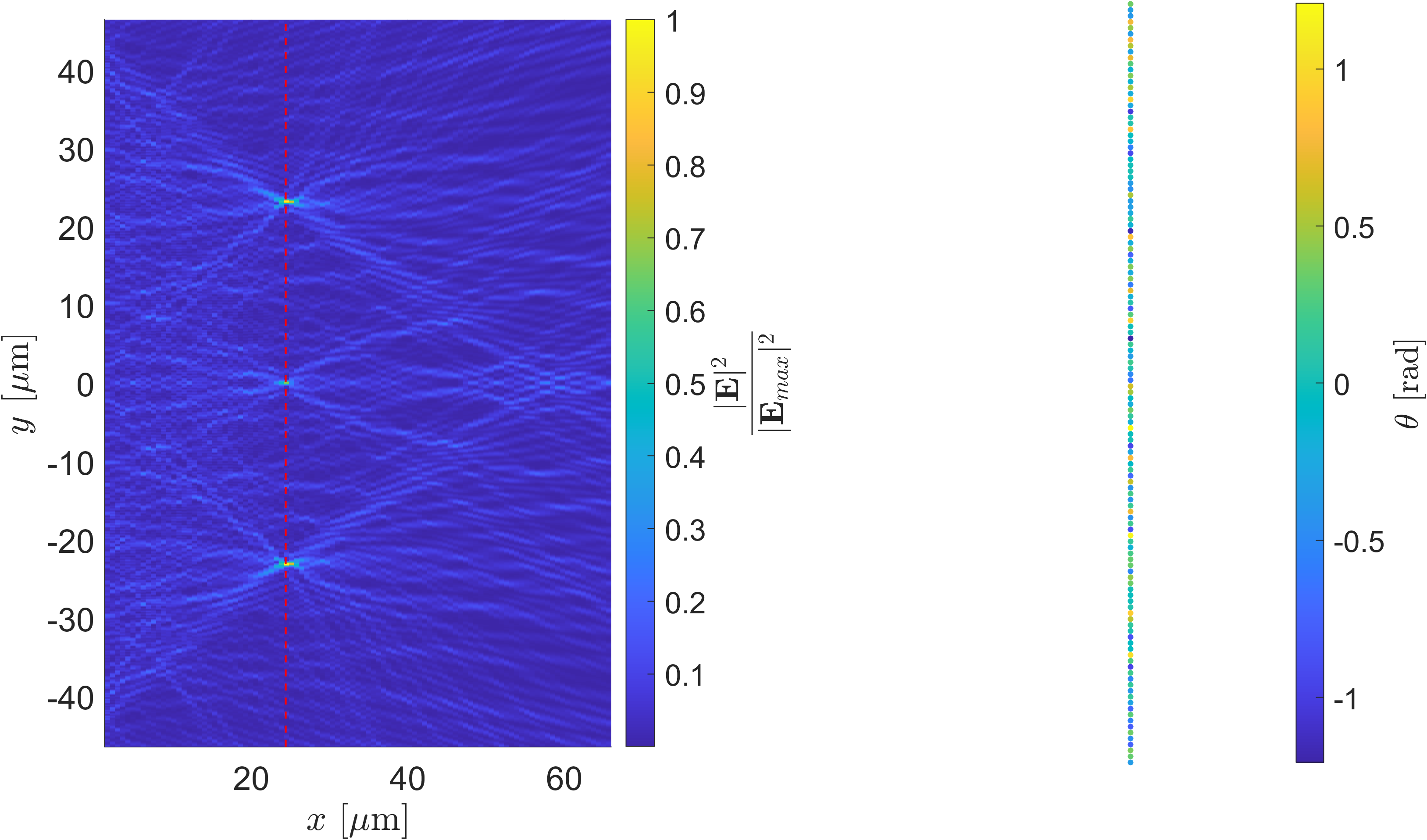}
    \caption{Light focusing in two focal spots. Left: Normalized intensity of the total electric field. The red dashed line is at $37\lambda_0$. Right: Optimal values of the rotation angle for each scatter.}
    \label{fig_res3_1}
\end{figure}

Similarly, Figure \ref{fig_res3_1} reports the results for three prescribed focal points located at $(37\lambda_0,35.18\lambda_0)$, $(37\lambda_0,0)$ and $(37\lambda_0,-35.18\lambda_0)$. The corresponding target field is defined by extending \eqref{eqn_super_dip} to include the superposition of the fields radiated by three dipoles positioned at the desired focal points.

In all cases, the light is effectively concentrated at the intended locations through rapidly varying meta-surface patterns, demonstrating the capability to handle more complex multi-spot focusing tasks with satisfactory accuracy.

\section{Conclusions}\label{sec:conclusions}
The inverse design of large-scale meta-surfaces is a challenging task, as the computational cost associated with their electromagnetic analysis is extremely high. In this paper, the adjoint variable method is investigated as an efficient tool for the fast evaluation of gradients required in most iterative optimization schemes.

The adjoint variable method has been extensively applied to freeform optimization problems, where the boundary of each meta-atom can vary continuously to minimize a given cost function. Although promising results have been demonstrated on small-scale structures, the high computational burden associated with binarization limits the scalability of these approaches.

In this work, a different strategy is proposed. The initial shape of each meta-atom is fixed, and a relatively small set of geometric parameters (such as width, height, rotation angle, or position) is chosen as the optimization variables. This significantly reduces the number of unknowns in large-scale problems.

The proposed parameterization is then combined with the adjoint variable method, enabling the computation of the gradient with respect to all design parameters through the solution of only one direct problem. Unlike previous works, the gradient is derived with respect to a general affine transformation, meaning that each meta-atom undergoes anisotropic scaling, rotation, and translation at each iteration of the optimization process. This allows for complex and non-trivial designs to be achieved without increasing the computational cost of the optimization.

More specifically, it is shown that the gradient computation via the adjoint method can be decoupled into two independent problems: (i) the determination of a suitable weighting function accounting for the considered geometric transformation, and (ii) the identification of the appropriate source term for the adjoint problem. A general framework is presented for deriving the weighting function in arbitrary scenarios, and explicit expressions are provided for the case of affine transformations. In addition, the adjoint sources corresponding to several representative cost functions are derived and discussed.

Numerical results confirm the feasibility and accuracy of the proposed method, showing excellent agreement between the gradients computed via the adjoint variable method and those obtained using finite differences. Furthermore, inverse design examples are presented for light-focusing applications, demonstrating the potential of the proposed approach.

\section{Acknowledgment}
This work was supported by the Italian Ministry of University and Research under the PRIN-2022, Grant Number 2022Y53F3X \lq\lq Inverse Design of High-Performance Large-Scale Metalenses\rq\rq.

\section*{Authorship contribution statement}

{\bf A. Tamburrino}: Conceptualization, Methodology, Formal analysis, Writing, Supervision.

{\bf V. Mottola}: Conceptualization, Methodology, Formal analysis, Writing.

\clearpage

\bibliographystyle
{iopart-num}
\bibliography{bibliography}
\end{document}